\documentclass[a4paper,12pt,reqno]{amsart}

\usepackage{amscd}%
\usepackage{amsmath}%
\usepackage{amsfonts}%
\usepackage{amssymb}%
\theoremstyle{plain}
\newtheorem{theorem}{Theorem}[section]
\newtheorem{lemma}{Lemma}[section]

\newtheorem{definition}{Definition}[section]
\newtheorem{corollary}{Corollary}[section]

\title{Isometric representations of totally ordered semigroups}
\author{M.A.Aukhadiev}
\address[M.A.Aukhadiev, V.H.Tepoyan]
{Kazan State Power Engineering University, Krasnoselskaya str., 51, 420066, Kazan, Russia }
\email[M.A.Aukhadiev]{m.aukhadiev@gmail.com}%
\author{V.H.Tepoyan}
\email[V.H.Tepoyan]{tepoyan.math@gmail.com}

\keywords{totally ordered semigroup, group, inverse semigroup, regular representation, isometric representation.}%
\begin{document}

\begin{abstract}
Let S be a subsemigroup of an abelian torsion-free group G. If S is a positive cone of G, then all C*-algebras generated by faithful isometrical non-unitary representations of S are canonically isomorphic. Proved by Murphy, this statement generalized the well-known theorems of Coburn and Douglas. In this note we prove the reverse. If all C*-algebras generated by faithful isometrical non-unitary representations of S are canonically isomorphic, then S is a positive cone of G. Also we consider $G=\mathbb{Z}\times\mathbb{Z}$ and prove that if S induces total order on G, then there exist at least two unitarily not equivalent irreducible isometrical representation of S. And if the order is lexicographical-product order, then all such representations are unitarily equivalent. 
\end{abstract}
\maketitle

\section{Introduction and preliminaries}
  
Within this paper $S$ is a subsemigroup of an additive abelian torsion-free group $G$ with zero. $S$ induces a partial order on $G$: $a\prec b$ if there exists $c\in S$ such that $a+c=b$. Semigroup $S$ induces full order on $G$, i.e. for any $a,b\in S$ either $a\prec b$ or  $b\prec a$, if $G=S\cup (-S)$ and $S\cap (-S)=\{0\}$. In this case write $S=G^+$ -- \emph{a positive cone} of $G$. Each semigroup $S$, which doesn't contain groups, is contained in some positive cone $G^+$. This follows from the axiom of choice.

Let $G$ be an abelian totally ordered group and $S$ -- subsemigroup of $G^+$, which doesn't contain groups. We denote by $\Delta_S$ a set of unitary equivalence classes of faithful irreducible non-unitary isometrical representations of semigroup $S$. For $V\in\Delta_S$ define $S_V$ as a semigroup generated by operators $V_{a}$ and $V^{*}_{b}$, where $a,b\in{S}$ and  $V_a=V(a)$.

An \emph{inverse semigroup} $P$ is a semigroup, such that each element $x$ has a unique \emph{inverse} element $x^*$, which satisfies the following:
\begin{equation}\label{inv}
xx^*x=x, \ x^*xx^*=x^*
\end{equation}

\begin{definition}
We call the representation $V\in{\Delta_S}$ inverse, if $S_V$ is an inverse semigroup.
\end{definition}

In the well-known work \cite{Coburn} Coburn proved that all isometric representations of semigroup $\mathbb{N}$ generate canonically isomorphic $C^*$-algebras. The same was proved by Douglas \cite{Douglas} for positive cones in $\mathbb{R}$ and by Murphy \cite{Murphy} for positive cones of abelian totally ordered groups. In section 2 we show that every semigroup $S$ has at least one inverse representation. Therefore all faithful isometric representations of  positive cone are inverse.

S.A.Grigoryan assumed that all representations in  $\Delta_S$ are inverse if and only if $S$ is a totally ordered semigroup, i.e. $S$ is a positive cone of some group. We prove this hypothesis in section 2.

In section 3 we prove that if $S$ induces full archimedian order on $\mathbb{Z}\times\mathbb{Z}$, then it has at least two unitarily not equivalent irreducible isometric representations. In case $S$ induces a total lexicographical-product order, all such representations are unitarily equivalent.

\section{Inverse representations}

\emph{Regular isometric representation} is a map $V:S\rightarrow{B(l^2(S))}$, $a\mapsto V_a$, defined as follows:

\[
(V_af)(b) =
\begin{cases}
f(c), & \text{if $b=a+c$ for some $c\in S$;} \\
0, & \text{ otherwise }
\end{cases}
\]

$C^*$-algebra generated by regular isometric representation of semigroup $S$ is called a \emph{reduced semigroup $C^*$-algebra}, denoted by $C^{*}_{red}(S)$ \cite{Jang}.

A finite product of operators of the form $V_{a}$ and $V^{*}_{b}$, $a,b\in{S}$ is called a \emph{monomial}. \emph{An index of monomial} $W=V_{a_1}V^{*}_{a_2}V_{a_3}...V^{*}_{a_n}$ is an element of group $\Gamma=S-S$, equal to $$\mathrm{ind} W= (a_2+a_4+...+a_n)-(a_1+a_3+...+a_{n-1}),$$ when $n$ is even \cite{Salahutdinov}. For odd $n$ we have:
 $$W=V_{a_1}V^{*}_{a_2}V_{a_3}...V_{a_n},$$ $$\mathrm{ind}W=(a_2+a_4+...+a_{n-1})-(a_1+a_3+...+a_{n}).$$ It is clear that $$\mathrm{ind}(W_1\cdot{W_2})=\mathrm{ind}W_1+\mathrm{ind}W_2.$$

Due to definition, monomials form a semigroup, which we denote by $S_V$.

\begin{lemma}\label{1}
The regular isometric representation of $S$ is inverse.
\begin{proof}
Consider a family $\left\{e_a\right\}_{a\in{S}}$ of elements in $l^2(S)$ such that $e_a(b)=\delta_{a,b}$. This is a natural orthonormal basis in $l^2(S)$. Every monomial $W$ in $S_V$ satisfies the following:
$$We_b=e_{b-d} \ \mbox{ or } \ 0, \ \mbox{ where } \ d=\mathrm{ind}W.$$
  
Note that $WW^*$ and $W^*W$ are monomials also, besides $$\mathrm{ind}(W\cdot W^*)=\mathrm{ind}(W^*\cdot W).$$
 By virtue of Lemma~2.2 in \cite{Salahutdinov2}, $WW^*$ and $W^*W$ are orthogonal projections. This implies immediately that $W=WW^*W$ and $W^*=W^*WW^*$. Therefore, an inverse element for $W$ is $W^*$.

\end{proof}

\end{lemma}

\begin{lemma}\label{lemma2}
There exists at least one noninverse representation in $\Delta_S$ for a semigroup $S\subsetneq{G^+}$.
\begin{proof}
Take a regular representation $V$ of $S$ in $B(l^2(S))$, $a\mapsto V_a$. Since $S$ is not equal to $G^+$, there exist incomparable elements $c,d\in S$, i.e. $c-d\notin S$ and $d-c\notin S$. Consider function $g_{c,d}=\frac{e_c+e_d}{\sqrt{2}}$ in $l^2(S)$. Denote by $H$ a Hilbert space generated by linear span of $\{V_ag_{c,d}\}_{a\in S}$. Note that $V_ag_{c,d}=g_{c+a,d+a}$.
Define representation $\widetilde{V}$ of semigroup $S$ on $H$, $a\mapsto \widetilde{V_a}$, by setting $\widetilde{V_a}=V_aP$, where $P:l^2(S)\rightarrow H$ is a projection on $H$.

This representation is faithful isometric due to its definition.

Let us show that
\begin{equation}\label{eq11}
\widetilde{V^{}_{c}}\widetilde{V^{*}_{c}}\widetilde{V^{}_{d}}\widetilde{V^{*}_{d}}\neq\widetilde{V^{}_{d}}\widetilde{V^{*}_{d}}\widetilde{V^{}_{c}}\widetilde{V^{*}_{c}}.
\end{equation}
Consider $\widetilde{V^{*}_{d}}g_{2c,c+d}$ and find such elements $x\in S$ that 
$$(\widetilde{V^{*}_{d}}g_{2c,c+d},g_{c+a,d+a})=0.$$

To this end, calculate

\begin{equation}\label{eq22}\begin{array}{c}
(\widetilde{V^{*}_{d}}g_{2c,c+d},g_{c+a,d+a})=(g_{2c,c+d},g_{c+d+a,2d+a})=\\  
 =(\frac{e_{2c}+e_{c+d}}{\sqrt{2}},\frac{e_{c+d+a}+e_{2d+a}}{\sqrt{2}})=\\
=\frac{1}{2}((e_{2c},e_{c+d+a})+(e_{2c},e_{2d+a})+(e_{c+d},e_{c+d+a})+(e_{c+d},e_{2d+a})).
\end{array}\end{equation}

%

First and last summands are equal to zero, since $c$ and $d$ are incomparabe. Therefore the scale product $(\widetilde{V^{*}_{d}}g_{2c,c+d},g_{c+a,d+a})$ is not equal to zero if and only if either $a=0$ or $a=2c-2d$. Note that element $2c-2d$ may not be contained in semigroup $S$. Despite this fact we continue the proof assuming $2c-2d\in S$. One can easily see that without this assumption the proof is trivial.

Denote by $H_0$ a Hilbert space in $H$ generated by elements of the following set $$\{g_{c+a,d+a} | \ a\neq 0, a\neq 2c-2d\}$$ Repeating the same arguments as above one can show that  $g_{c,d}$ and $g_{3c-d,2c-d}$ are mutually orthogonal, and both are orthogonal to $H_0$. Consequently, $codimH_0=2$ and the elements $g_{c,d}$ and $g_{3c-d,2c-d}$ form an orthonormal basis in $H^{\bot}_{0}\subset H$. Thus,
$$H=H_0\oplus\mathbb{C}g_{c,d}\oplus\mathbb{C}g_{3c-d,2c-d},$$
and from equation (\ref{eq22}) we have
$$V^{*}_{d}g_{2c,c+d}=\frac{1}{2}(g_{c,d}+g_{3c-2d,2c-d}).$$

For futher be noted, the assumption $2c-2d\in S$ implies that $2d-2c$ is not contained in semigroup $S$. Otherwise $G^+$ would contain non-trivial group, which is impossible. Therefore, due to symmetry we get $$V^{*}_{c}g_{c+d,2d}=\frac{1}{2}g_{c,d}.$$

Thus, 
$$\widetilde{V^{}_{c}}\widetilde{V^{*}_{c}}\widetilde{V^{}_{d}}\widetilde{V^{*}_{d}}g_{2c,c+d}=\frac{1}{2}\widetilde{V^{}_{c}}\widetilde{V^{*}_{c}}\widetilde{V^{}_{d}}(g_{c,d}+g_{3c-2d,2c-d})=$$ $$=\frac{1}{2}(\widetilde{V^{}_{c}}\widetilde{V^{*}_{c}}g_{c+d,2d}+\widetilde{V^{}_{c}}\widetilde{V^{*}_{c}}g_{3c-d,2c})=$$
$$=\frac{1}{4}V_cg_{c,d}+\frac{1}{2}V_cg_{2c-d,c}=\frac{1}{4}g_{2c,c+d}+\frac{1}{2}g_{3c-d,2c}.$$

On the other hand,
$$\widetilde{V_{d}}\widetilde{V^{*}_{d}}\widetilde{V_{c}}\widetilde{V^{*}_{c}}g_{2c,c+d}=\widetilde{V_{d}}\widetilde{V^{*}_{d}}g_{2c,c+d}=\frac{1}{2}g_{c+d,2d}+\frac{1}{2}g_{3c-d,2c}.$$

Consequently, we get inequality (\ref{eq11})

\end{proof}

\end{lemma}

\begin{theorem}
The following properties of semigroup $S$ are equivalent
\begin{enumerate}
	\item $S=G^+$;\label{a}
	\item all representations in $\Delta_S$ are canonically isomorphic;\label{b}
	\item all representations in $\Delta_S$ are inverse;\label{c}
	\item for any representation $V$ in $\Delta_S$ and for any $a,b\in S$ the following equality is satisfied\label{d}
	$$V^{}_{a}V^{*}_{a}V^{}_{b}V^{*}_{b}=V^{}_{b}V^{*}_{b}V^{}_{a}V^{*}_{a}.$$
\end{enumerate}

\begin{proof}
(\ref{a})$\Rightarrow$(\ref{b}) was proved by Murphy \cite{Murphy}.

Let us show implication (\ref{b})$\Rightarrow$(\ref{c}). Suppose all representations in $\Delta_S$ are canonically isomorphic and $S\subset{G^+}$. Consider representation $V:S\rightarrow{l^2(G^+)}$, $a\mapsto{V_a}$, defined by
$$V_ae_b=e_{a+b},$$
where $\left\{e_a\right\}_{a\in{G^+}}$ is an orthonormal basis in $l^2(G^+)$. For any $a,b,c\in S$ if $a\prec b$ or $a=b$ we have $V^{*}_{a}V^{}_{b}e^{}_{c}=e^{}_{c+b-a}$.  Since all elements in $G^+$ are pairwise comparable, we have two cases. If $a\prec{b}$, then operator $V^{*}_{a}V^{}_{b}$ is isometric, otherwise ($b\prec a$) operator $(V^{*}_{a}V^{}_{b})^*$ is isometric. Consequently, semigroup $S_V$ is inverse.

Implication (\ref{c})$\Rightarrow$(\ref{d}) concerns only inverse semigroups, and it was proved in \cite{Clifford}.

Lemma ~\ref{lemma2} implies (\ref{d})$\Rightarrow$(\ref{a}).

\end{proof}

\end{theorem}

\begin{corollary}\label{cor1}
The $C^*$-algebras $C^{*}(S)$ and $C^{*}_{red}(S)$ are isomorphic if and only if $S$ is totally ordered, where $C^{*}_{}(S)$ is a universal enveloping $C^{*}$-algebra, generated by all isometric representations of semigroup $S$ \cite{Murphy2}. 
\end{corollary}

In particular case $S=\mathbb{Z}^+$ this statement implies that the algebras $C^{*}_{}(\mathbb{Z}^+)$ and $C^{*}_{red}(\mathbb{Z}^+)$ are isomorphic. This result was proved by Coburn in his well-known work \cite{Coburn}. As an example of the converse to this statement take $S=\mathbb{Z}^+\backslash \left\{1\right\}$. Due to Corollary~\ref{cor1}, the algebras $C^{*}(S)$ and $C^{*}_{red}(S)$ are not isomorphic. The same was shown in \cite{Jang2}, and  this case was studied in details in \cite{Raeburn}.

\section{$C^*$-algebras generated by totally ordered semigroup in $\mathbb{Z}\times\mathbb{Z}$}

Consider group $G=\mathbb{Z}\times\mathbb{Z}$. Total order on $G$ is equivalent to straight line dividing it into two parts. It implies two cases. The first case: the line meets the point $(0,0)$ and doesn't meet any integers. Such line is characterized by equation $x+\alpha y=0$, where $\alpha$ is irrational. The second case: the line meets integers, i.e. $x+\alpha y=0$, for rational $\alpha$. The order induced by the first line is archimedian. In the second case we may consider $G=S\cup(-S)$, where $S=\left\{(n,m)\in\mathbb{Z}\times\mathbb{Z} \  | \ m>0 \  \mbox{ or } (n,0), \  n\geq 0\right\}$ $(S=G^+)$. In this case the order cannot be archimedian, since we have $(-1,1)<(0,1)$ together with $n\cdot(-1,1)<(0,1)$ for any $n>0$.

\begin{theorem}
 \ \  
\begin{enumerate}
	\item If $G^+$ induces total archimedian order on $G$, then $card\Delta_S>1$;\label{i}
	\item If $G^+$ induces lexicographical-product order, then $card\Delta_S=1$.\label{i2}
\end{enumerate}

\begin{proof}
(\ref{i}) Suppose $G^+$ induces total archimedian order on $G$. Without loss of generality, we may assume that $G^+\subset{\mathbb{R}^+}$. Therefore  $\overline{G^+}=\mathbb{R}^+$. Let us give a new representation of semigroup $G^+$.

Consider the Hardi space $H^2$. By the help of inner singular function $exp\{\frac{1+e^{i\theta}}{1-e^{i\theta}}\}$ define nonunitary faithful isometric representation of the semigroup $\mathbb{R}^+$ in $B(H^2)$, $t\mapsto V_t$, by the following equation:
$$(V_tg)(e^{i\theta})=exp(t\frac{1+e^{i\theta}}{1-e^{i\theta}})g(e^{i\theta}).$$

One can easily verify that $V_t$ is an isometric operator on $H^2$. Let us show that this representation is not uquivalent to regular representation.

In case of regular representation $W$ there exists element $e_0$ such that $W_te_0\bot{e_0}$ for any $t\in{G^+}$. It is sufficient to show that $H^2$ does not contain element $g$, such that $V_tg\bot{g}$ for any $t\in{G^+}$. Indeed, suppose that there exists such element $g$. Then we have

\begin{equation}\label{int1}
0=(V_tg,g)=\frac{1}{2\pi}\int_{S^1}exp(t\frac{1+e^{i\theta}}{1-e^{i\theta}})g(e^{i\theta})\overline{g(e^{i\theta})}d\mu(\theta)= \end{equation}
$$=\frac{1}{2\pi}\int_{S^1}exp(t\frac{1+e^{i\theta}}{1-e^{i\theta}})d\mu(\theta).$$

If $t\rightarrow 0$, the right-hand side of (\ref{int1}) converges to 1, which leads to a contradiction. Thus, representations $V$ and $W$ are not equivalent. 

Now let us prove the second part of the theorem, (\ref{i2}).

The group of transformations of iteger lattice $\mathbb{Z}\times\mathbb{Z}$ is a group $SL(2,\mathbb{Z})$. For any pair of lexicographical-product orders on $\mathbb{Z}\times\mathbb{Z}$ there exists an element of $SL(2,\mathbb{Z})$, which transforms the first one to the second one. Therefore, without loss of generality, we may consider that $S$ is equal to the following semigroup:
$$\{(n,m)\in\mathbb{Z}\times\mathbb{Z} \  | \ m>0 \  \mbox{ or } (n,0), \  n\geq 0\}.$$

Take representation $V:S\rightarrow B(H)$ in $\Delta_S$. Since operator $V_{(1,0)}$ is isometric and not unitary, there exists $h_0\in H$ such that $V^{*}_{(1,0)}h_0=0$. Since $V_{(0,1)}=V_{(1,0)}V_{(-1,1)}$, we have 
\begin{equation}\label{eq32}
V^{*}_{(0,1)}h^{}_{0}=V^{*}_{(-1,1)}V^{*}_{(1,0)}h^{}_{0}=0.
\end{equation}
Therefore, $h_0$ is an initial vector for operators $V_{(0,1)}$ and $V_{(1,0)}$. Consequently, it is initial for any $V_{(n,m)}$, where $(n,m)\in S$.

Consider Hilbert space $H_1$, generated by linear span of the set $$\{V_{(n,m)}h_0, \ (n,m)\in S\}$$ Equation (\ref{eq32}) implies that the family $\{V_{(n,m)}h_0, \ (n,m)\in S\}$ forms an orthonormal basis in $H_1$, and
$$V^{*}_{(k,l)}V^{}_{(n,m)}h_0=V^{}_{(a,b)}h_0 \ \mbox{ or } 0.$$
Therefore, $H_1$ is an invariant subspace for $C^*$-algebra $C^{*}_{red}(S)$. Since representation $V$ is irreducible, we have $H_1=H$.

Consequently, the family of vectors $e_{n,m}=V_{(n,m)}h_0$, for $(n,m)\in S$, forms an orthonormal basis of $H$. This implies immediately $H\cong l^2(S)$.

\end{proof}

\end{theorem}

\end{document}